\documentclass[12pt,oneside]{article}
\usepackage{amsmath,amssymb,amsfonts,amsthm}
\usepackage{color}

\newcommand{\T}{{\cal T}}

\newcommand{\Real}{\mathbb R}

\newcommand{\To}{\longrightarrow}

\newcommand{\prof}{\noindent \textit{\textbf{Proof.\:\:}}}
\newcommand{\tm}{\T M}
\newcommand{\p}{\pi^{-1}(TM)}
\newcommand {\cp}{\mathfrak{X}(\pi (M))}
\newcommand {\ccp}{\mathfrak{X}^{*}(\pi (M))}

\def\o#1{\overline{#1}}

\def\f{ \mathfrak{F}(M)}

\def\Section#1{\vspace{30truept}\addtocounter{section}{1}\setcounter{thm}{0}\setcounter{equation}{0}
{\noindent\Large\bf\arabic{section}.~~#1}\par \vspace{12pt}}

\newtheorem{thm}{Theorem}[section]
\newtheorem{cor}[thm]{Corollary}
\newtheorem{lem}[thm]{Lemma}
\newtheorem{prop}[thm]{Proposition}
\newtheorem{defn}[thm]{Definition}

\numberwithin{equation}{section}
\begin{document}
%\title{{\bf CHARACTERIZATION OF CLOSED VECTOR FIELDS IN FINSLER GEOMETRY}}

\title{{\bf Characterization of Closed Vector Fields in Finsler Geometry}}
\author{\bf{Nabil L. Youssef }}
\date{}
\maketitle

% End of preamble and beginning of text.
% Produces the title.
\vspace{-1.15cm}

\begin{center}
{Department of Mathematics, Faculty of Science,\\
Cairo University, Giza, Egypt.\\
nyoussef@frcu.eun.eg,\, nlyoussef2003@yahoo.fr}
\end{center}
%\smallskip
\begin{center}
{\bf Dedicated to the memory of Prof. Dr. A. Tamim}
\end{center}

\bigskip

\noindent{\bf Abstract.}  The $\pi$-exterior derivative $\o d$,
which is the Finslerian generalization of the (usual) exterior
derivative $d$ of Riemannian geometry, is defined. The notion of a
$\o d$-closed vector field is introduced and investigated. Various
characterizations of \linebreak $\o d$-closed vector fields are
established. Some results concerning $\o d$-closed vector fields in
relation to certain special Finsler spaces are obtained.
\footnote{This paper was presented in \lq \lq \,The
 International Conference on Finsler Extensions of Relativity Theory \,\rq\rq\ held at
Cairo, Egypt, November 4-10, 2006.}



%%%%%%%%%%%%%%%%%%%%%%%%%%%%%%%%%%%%%%%%%%%%%% INTRODUCTION %%%%%%%%%%%%%%%%%%%%%%%%%%%%%%%%%%%%%%%%%%%%%%%%%%%%%%%%%%%%%%%

\vspace{30truept}\centerline{\Large\bf{Introduction}}\vspace{12pt}
\par
In the present work, we adopt the pullback approach to Finsler
geometry. In Finsler geometry, there is a canonical linear
connection (corresponding to the Levi-Civita connection of
Riemannian geometry), due to E. Cartan, which is not a\linebreak
connection on the manifold $M$ but is a connection on
$\,\pi^{-1}(TM) $, the pullback of the tangent bundle $TM$ by
$\,\pi: \T M\longrightarrow M$. The Cartan connection plays a key
role in this work.

\vspace{4pt}
\par We define the notion of $\pi$-exterior derivative $\o d$, which is a natural generalization to
Finsler geometry of the (usual) exterior derivative $d$ of
Riemannian geometry. We then introduce and investigate an important
class of $\pi$-vector fields on a Finslr manifold, which we refer to
as $\o d$-closed vector fields. Various characterizations of such
$\pi$-vector fields are established. Some results concerning $\o
d$-closed vector fields in relation to certain special Finsler
spaces are obtained. The notion of a $\pi$-distribution is also
introduced and is related to $\o d$-closed vector fields. It should
finally be noted that our investigation is entirely global or
intrinsic (free from local coordinates).

\vspace{4pt} The idea of this work is due to Prof. A. Tamim, whom we
miss profoundly.

%%%%%%%%%%%%%%%%%%%%%%%%%%%%%%%%%%%%%%%%% SECTION: PRELIMINARIES %%%%%%%%%%%%%%%%%%%%%%%%%%%%%%%%%%%%%%%%%%%%%%%%%%%%%%%%%%%%%%%%%%
\newpage

\Section{Notation and Preliminaries}

In this section, we give a brief account of the basic concepts
 of the pullback\linebreak formalism of Finsler geometry necessary for this work. For more details
  refer to~\cite{r57},\,\cite{r61} and~\,\cite{r44}.
 We make the general
assumption that all geometric objects we consider are of class
$C^{\infty}$.
\par
The following notations will be used throughout the present paper:\\
 $M$: a real differentiable manifold of finite dimension $n$ and of
class $C^{\infty}$,\\
 $\mathfrak{F}(M)$: the $\Real$-algebra of differentiable functions
on $M$,\\
 $\mathfrak{X}(M)$: the $\mathfrak{F}(M)$-module of vector fields
on $M$,\\
$\pi_{M}:TM\longrightarrow M$: the tangent bundle of $M$,\\
$\pi: \T M\longrightarrow M$: the subbundle of nonzero vectors
tangent to $M$,\\
$V(TM)$: the vertical subbundle of the bundle $TTM$,\\
 $P:\pi^{-1}(TM)\longrightarrow \T M$ : the pullback of the
tangent bundle $TM$ by $\pi$,\\
$P^*:\pi^{-1}(T^{*}M)\longrightarrow \T M$ : the pullback of the
cotangent bundle $T^{*}M$ by $\pi$,\\
 $\mathfrak{X}(\pi (M))$: the $\mathfrak{F}(TM)$-module of
differentiable sections of  $\pi^{-1}(T M)$,\\
$\mathfrak{X}^{*}(\pi (M))$: the $\mathfrak{F}(TM)$-module of
differentiable sections of  $\pi^{-1}(T^{*} M)$,\\
$ i_{X}$ : interior product with respect to  $X
\in\mathfrak{X}(M)$,\\
$df$ : the exterior derivative  of $f\in \mathfrak{F}(M)$,\\
$ d_{L}:=[i_{L},d]$, $i_{L}$ being the interior derivative with
respect to the vector form $L$. \vspace{8pt}
\par
Elements  of  $\mathfrak{X}(\pi (M))$ will be called $\pi$-vector
fields and will be denoted by barred letters $\overline{X} $. Tensor
fields on $\pi^{-1}(TM)$ will be called $\pi$-tensor fields. The
fundamental $\pi$-vector field is the $\pi$-vector field
$\overline{\eta}$ defined by $\overline{\eta}(u)=(u,u)$ for all
$u\in \T M$. The lift to $\pi^{-1}(TM)$ of a vector field $X$ on $M$
is the $\pi$-vector field $\overline{X}$ defined by
$\overline{X}(u)=(u,X(\pi (u)))$. The lift to $\pi^{-1}(TM)$ of a
$1$-form $\omega$ on $M$ is the $\pi$-form $\overline{\omega}$
defined by $\overline{\omega}(u)=(u,\omega(\pi (u)))$. \vspace{8pt}
\par
We have the following short exact sequence of vector bundles,
relating the tangent bundle $T(\T M)$ and the pullback bundle
$\pi^{-1}(TM)$:\vspace{-0.2cm}
$$0\longrightarrow
 \pi^{-1}(TM)\stackrel{\gamma}\longrightarrow T(\T M)\stackrel{\rho}\longrightarrow
\pi^{-1}(TM)\longrightarrow 0 ,\vspace{-0.2cm}$$
 where the bundle morphisms $\rho$ and $\gamma$ are defined respectively by
$\rho = (\pi_{\T M},d\pi)$ and $\gamma (u,v)=j_{u}(v)$, where
$j_{u}$  is the natural isomorphism $j_{u}:T_{\pi_{M}(v)}M
\longrightarrow T_{u}(T_{\pi_{M}(v)}M)$.\linebreak The vector
$1$-form $J$ on $TM$ defined by $J=\gamma o \rho\,$ is called the
natural almost tangent structure of $T M$ and the vertical vector
field $\mathcal{C}$ on $TM$ defined by $\mathcal{C}:=\gamma o
\overline{\eta} $ is called the fundamental or the canonical
(Liouville) vector field \,\cite{r27}. \vspace{8pt}
\par
 Let $\nabla$ be  a linear
connection (or simply a connection) in the pullback bundle
$\pi^{-1}(TM)$.
 We associate to
$\nabla$ the map\vspace{-0.2cm}
$$K:T \T M\longrightarrow \pi^{-1}(TM):X\longmapsto \nabla_X \overline{\eta}
,\vspace{-0.2cm}$$ called the connection (or the deflection) map of
$\nabla$. A tangent vector $X\in T_u (\T M)$ is said to be
horizontal if $K(X)=0$ . The vector space $H_u (\T M)= \{ X \in T_u
(\T M) : K(X)=0 \}$ of the horizontal vectors
 at $u \in  \T M$ is called the horizontal space to $M$ at $u$  .
   The connection $\nabla$ is said to be regular if
 $$T_u (\T M)=V_u (\T M)\oplus H_u (\T M) \qquad \forall u\in \T M.$$
  \par
  If $M$ is endowed with a regular connection, then the vector bundle
   maps
\begin{eqnarray*}
 \gamma &:& \pi^{-1}(T M)  \To V(\T M), \\
   \rho |_{H(\T M)}&:&H(\T M) \To \pi^{-1}(TM), \\
   K |_{V(\T M)}&:&V(\T M) \To \pi^{-1}(T M)
\end{eqnarray*}
 are vector bundle isomorphisms.
   Let us denote
 $\beta=(\rho |_{H(\T M)})^{-1}$,
then \vspace{-0.2cm}
   \begin{align}\label{fh1}
    \rho  o  \beta = id_{\pi^{-1} (\T M)}, \quad  \quad
       \beta o \rho =\left\{
                                \begin{array}{ll}
                                          id_{H(\T M)} & \hbox{on   H(\T M)} \\
                                         0 & \hbox{on    V(\T M)}
                                       \end{array}
                                     \right.\vspace{-0.2cm}
\end{align}
\par
The classical  torsion tensor $\textbf{T}$  of a regular connection
$\nabla$ is defined by
$$\textbf{T}(X,Y)=\nabla_X \rho Y-\nabla_Y\rho X -\rho [X,Y] \quad
\forall\,X,Y\in \mathfrak{X} (\T M).$$ The horizontal ((h)h-) and
mixed ((h)hv-) torsion tensors, denoted by $Q $ and $ T $\linebreak
respectively, are defined by
$$Q (\overline{X},\overline{Y})=\textbf{T}(\beta \overline{X}\beta \overline{Y}),
\, \,\,\,\, T(\overline{X},\overline{Y})=\textbf{T}(\gamma
\overline{X},\beta \overline{Y}) \quad \forall \,
\overline{X},\overline{Y}\in\mathfrak{X} (\pi (M)).$$
\par The classical curvature tensor  $\textbf{K}$ of the connection
$\nabla$ is defined by
 $$ \textbf{K}(X,Y)\rho Z=-\nabla_X \nabla_Y \rho Z+\nabla_Y \nabla_X \rho Z+\nabla_{[X,Y]}\rho Z
  \quad \forall\, X,Y, Z \in \mathfrak{X} (\T M).$$
The horizontal (h-), mixed (hv-) and vertical (v-) curvature
tensors, denoted by $R$, $P$ and $S$ respectively, are defined by
$$R(\overline{X},\overline{Y})\o Z=\textbf{K}(\beta
\overline{X}\beta \overline{Y})\o Z,\,\,\,\,
P(\overline{X},\overline{Y})\o Z=\textbf{K}(\beta
\overline{X},\gamma \overline{Y})\o Z,\,\,\,\,
S(\overline{X},\overline{Y})\o Z=\textbf{K}(\gamma
\overline{X},\gamma \overline{Y})\o Z.$$ The contracted curvature
tensors, denoted by $\widehat{R}$, $\widehat{P}$ and $\widehat{S}$
respectively, are also known as the
 (v)h-, (v)hv- and (v)v-torsion tensors and are defined by
$$\widehat{R}(\overline{X},\overline{Y})={R}(\overline{X},\overline{Y})\o \eta,\quad
\widehat{P}(\overline{X},\overline{Y})={P}(\overline{X},\overline{Y})\o
\eta,\quad
\widehat{S}(\overline{X},\overline{Y})={S}(\overline{X},\overline{Y})\o
\eta.$$

Now, let $(M,L)$ be a Finsler manifold, where $L$ is the Lagrangian
defining the Finsler structure on $M$. Let $g$ be the Finsler metric
in $\p$  defined by $L$.

\begin{thm} \label{th.1}{\em{\cite{r44}}} Let $(M,L)$ be a Finsler manifold. There exists a
unique regular connection $\nabla$ in $\pi^{-1}(TM)$ such that
\begin{description}
  \item[(a)]  $\nabla$ is  metric\,{\em:} $\nabla g=0$,

  \item[(b)]   The horizontal torsion of $\nabla$ vanishes\,{\em:} $Q=0
  $,
  \item[(c)]  The mixed torsion $T$ of $\nabla$ satisfies \,
  $g(T(\overline{X},\overline{Y}), \overline{Z})=g(T(\overline{X},\overline{Z}),\overline{Y})$.
\end{description}
\end{thm}
\vspace{-3pt}
 Such a connection is called the Cartan connection
associated to the Finsler\linebreak manifold $(M,L)$. \vspace{6pt}
\par One can
show that the torsion of the Cartan connection has the property that
$T(\overline{X},\overline{\eta})=0$ for all $\overline{X} \in
\mathfrak{X} (\pi (M))$. For the Cartan connection, we have
\vspace{-0.2cm}
   \begin{align}\label{fh2}
    Ko\gamma = id_{\pi^{-1} (\T M)}, \quad \quad
      \gamma o K =\left\{
                                       \begin{array}{ll}
                                          id_{V(\T M)} & \hbox{on   V(\T M)} \\
                                         0 & \hbox{on    H(\T M)}
                                       \end{array}
                                     \right.
\end{align}
   Then, from (\ref{fh1}) and  (\ref{fh2}), we get\vspace{-0.2cm}
 \begin{equation}\label{hv}
   \beta o \rho + \gamma o K = id_{T(\T M)}\vspace{-0.2cm}
\end{equation}
 Hence, if we set
   $h {\!\!}:= \beta o \rho $\, and \,$v {\!\!}:=\gamma o K$, then every vector field $X{\!\!}\in{\!\!}\mathfrak{X}(\T M)$ can be
   represented uniquely in the form\vspace{-0.2cm}
           \begin{equation}\label{uy}
                                         X = hX +vX =\beta \rho X + \gamma K X \vspace{-0.2cm}
                                  \end{equation}
 The maps $h$ and $v$ are the horizontal and  vertical projectors associated
 with the Cartan connection $\nabla$: $h^{2}= h$, \,\,\,$v^{2}= v$,\,\,\,
  $h+v= id_{\mathfrak{X}(TM)}$,\,\, \, $voh=hov=0 .$

\begin{defn} {\em{\cite{r48}}} With respect to the Cartan  connection $\nabla$, we have:\newline
-- The horizontal Ricci tensor $Ric^h$ is defined by
 $$Ric^h(\overline{X},\overline{Y}):= Tr\{ \overline{Z} \longmapsto
R(\overline{X},\overline{Z})\overline{Y}\}, \quad \text {for all}
\quad \overline{X},\overline{Y}\in \mathfrak{X} (\T M).$$ -- The
horizontal Ricci  map $Ric_0^h$ is defined by
$$g(Ric_0^h(\overline{X}),\overline{Y}):=Ric^h(\overline{X},\overline{Y}),
\quad \text {for all} \quad \overline{X},\overline{Y}\in
\mathfrak{X} (\T M).$$
    --The horizontal scalar
curvature $Sc^h$ is defined by\vspace{-0.2cm}
     $$Sc^h:= Tr( Ric_0^h).$$
 \end{defn}

\par We terminate this section
by some concepts and results concerning the notion of a nonlinear
connection in the sense of Klein-Grifone\,\,\cite{Grifone
1},\,\,\cite{r27}: \vspace{-0.2cm}
\begin{defn}  A nonlinear connection on $M$ is a vector $1$-form $\Gamma$
on $TM$, $C^{\infty}$ on $\T M$, $C^{o}$ on $TM$, such
that\vspace{-0.3cm}$$J \Gamma=J, \quad\quad \Gamma J=-J .$$
\end{defn}
\vspace{-0,2cm}
 The  horizontal and vertical projectors $h$  and $v$ associated with $\Gamma$ are defined by
   $h:=\frac{1}{2} (I+\Gamma),\, v:=\frac{1}{2}
 (I-\Gamma).$
Thus $\Gamma$ gives rise to the decomposition $TTM= H(TM)\oplus
V(TM)$, where $H(TM):=Im \, h = Ker\, v $, $V(TM):= Im \, v=Ker \,
h$. The torsion $T$ of a nonlinear connection $\Gamma$ is the vector
$2$-form  on $TM$ defined by $T:=\frac{1}{2} [J,\Gamma]$.   The
curvature of a nonlinear connection $\Gamma$ is the vector $2$-form
$\Re$ on $TM$
     defined by
    $\Re:=- \frac{1}{2} [h,h].$
\begin{prop} Let $(M,L)$ be a Finsler manifold. The vector field $G$
 on $TM$ determined by $i_{G}\Omega =-dE$ is a spray,
 called the canonical spray associated with the energy $E$, where
 $E:=\frac{1}{2}L^{2}$ and $\Omega:=dd_{J}E$.
 \end{prop}
\vspace{-0.2cm}

\begin{thm} \label{th.9a} On a Finsler manifold $(M,L)$, there exists a unique
conservative nonlinear  connection {\em($d_{h}E=0$)} with zero
torsion. It is given by\,{\em:} \vspace{-0.3cm} $$\Gamma = [J,G]
,\vspace{-0.3cm} $$ where $G$ is the canonical spray.
\end{thm}
\vspace{-0.2cm} Such a connection is called the Barthel connection
or the Cartan nonlinear\linebreak connection associated with
$(M,L)$. \vspace{0.2cm}
\par It should be
noted that the horizontal and vertical projectors of the Cartan
connection and the barthel connection coincide. Also, the canonical
spray $G=\beta o \overline{\eta}$, and
 $G$ is thus horizontal.

%%%%%%%%%%%%%%%%%%%%%%%%%%%%%%%%%%%%%%%%%% SECTION: CLOSED VECTOR FIELDS %%%%%%%%%%%%%%%%%%%%%%%%%%%%%%%%%%%%%%%%%%%%%%%%%%%%%%%%%%%%%%%%%%%%%%%%
\newpage
\Section{$\pi$-Exterior derivative and $\o d$-closed vector field}
In this section, we introduce the notion of $\pi$-exterior
differentiation operator $\o d$, which is the Finslerian version of
the (usual) exterior differentiation operator $d$ of Riemannian
geometry. We then investigate the $\o d$-closed $\pi$-vector fields
in Finsler geometry.
\par Let $(M,L)$ be a Finsler manifold. Let $g$ be the Finsler
metric defined by the Lagrangian $L$ and $\nabla$ be the Cartan
connection associated with $(M,L)$.
\par
\vspace{0.2cm} We start with the following lemma which is useful for
subsequent use.

\begin{lem}\label{lem.bracket}For all
$\,\overline{X},\overline{Y}\in \mathfrak{X}(\pi(M))$ and
$X,Y\in\mathfrak{X}(TM)$,
 we have{\,\em:}\vspace{-0.2cm}

\begin{description}
    \item[(a)]$\widehat{R}(\o X, \o Y)=K[\beta \overline{X},\beta \overline{Y}],$

    \item[(b)]$ \mathfrak{R}(X,Y)=-\gamma \widehat{R}(\rho X, \rho Y),$
\end{description}
\vspace{-0.2cm} $\mathfrak{R}$ being the curvature tensor of the
Barthel connection.
\end{lem}

%\begin{proof}
\prof \\ \textbf{(a)} One can easily show that \vspace{-0.2cm}
\begin{equation}\label{2eq.5}
 [\beta \overline{X},\beta \overline{Y}]=
     \gamma(R(\overline{X},\overline{Y})\overline{\eta})
     + \beta(\nabla_{\beta \overline{X}}\overline{Y}-
     \nabla_{\beta \overline{Y}}\overline{X}).\vspace{-0.2cm}
\end{equation}
Then, (a) follows directly from (\ref{2eq.5}) by applying the
operator $K$ on both sides, noting that
$K\circ\gamma=id_{\pi^{-1}(TM)}$ and $K\circ\beta=0$. \\
\textbf{(b)} Setting $\o X= \rho X$ and $\o Y= \rho Y$ in (a), we
get\vspace{-0.2cm}
$$\gamma \hat{R}(\rho X, \rho Y)=\gamma K[\beta \rho{X},\beta \rho{Y}]=
v[hX,hY].\vspace{-0.2cm}$$ Then, (b) follows directly from the
identity \cite{r85}: \,\,$\mathfrak{R}(X,Y)=-v[hX,hY]$. \ \ $\Box$
%\end{proof}
\begin{defn}\label{2def.4}For any  given $\pi$-vector field $\o X$,
 the $\mathfrak{F}(TM)$-linear operator $A_{\o X}$ on
$\cp$ is defined by\,{\em:}\vspace{-0.1cm}
\begin{equation}\label{2eq.4}
   A_{\o X}(\o Y):=\nabla_{\beta \o Y}\o X.
\end{equation}
\end{defn}

\begin{defn}\label{2def.1} {\em \cite{r77}} Let $\omega$ be a $\pi$-form of order
$p\geq0$.\\ For $p>0$, we define\vspace{-0.2cm}
\begin{equation}\label{2eq.1}
          \begin{array}{rcl}
        (\o d \omega)(\o X_{1},...,\o X_{p+1})&: = &\sum_{i=1}^{p+1}(-1)^{i+1} \beta \o X_{i}
        \cdot \omega(\o X_{1},...,\widehat{\o X_{i}},...,\o X_{p+1}) \\
        & + & \sum_{i<j}(-1)^{i+j}\omega(\rho[\beta \o X_{i}, \beta \o X_{j}],\o X_{1},...,
        \widehat{\o X_{i}},...,\widehat{\o X_{j}},...,\o X_{p+1})\,\,\,\vspace{-0.2cm}
    \end{array}
  \end{equation}
 \noindent For  $p=0$, we set \vspace{-0.2cm}
 \begin{equation}\label{2eq.2}
    (\o d \omega)(\o X):=\beta \o X \cdot \omega, \ \text{that is},\
    \o d \omega= d \omega \,o \,\beta .\vspace{-0.2cm}
\end{equation}
{\em({\it Here, the symbol \lq\lq\, $\widehat{}$ \,\rq\rq \, means
that the corresponding argument is omitted.})}\\
In particular, for a {\em$(1)$}$\pi$-form $\omega$, we
have\vspace{-0.2cm}
\begin{equation}\label{2eq.3}
    (\o d \omega)(\o X, \o Y)=\beta \o X \cdot \omega(\o Y)-
    \beta \o Y \cdot \omega(\o X)- \omega(\rho[\beta \o X, \beta \o Y]) .\vspace{-0.2cm}
\end{equation}
The operator  $\o d$ will be called the $\pi$-exterior derivative.
\end{defn}

\begin{defn}\label{2def.2} A {\em$(p)$}$\pi$-form $\omega$ is said to
be: \\
-- $\o d$-closed if $\o d \omega=0$.\\
-- $\o d$-exact if $ \omega=\o d \alpha$, for some
{\em$(p-1)$}$\pi$-form $\alpha$.
\end{defn}
\par It should be noted that a $\o d$-exact $\pi$-form is not
necessarily  $\o d$-closed, contrary to the case of the (ordinary)
exterior derivative $d$ ({\cite{r80}} and {\cite{r82}}). This is due
to the fact that the property ${\o d}^{2}=0$ dose not hold for the
$\pi$-exterior derivative $\o d$.

\begin{defn}\label{2def.3}Let $\o X\in\cp$ be a $\pi$-vector field
and let $\omega\in\ccp$ be the associated {\em$(1)$}$\pi$-form under
the duality defined by the Finsler metric $g${\em\,:} $\omega=i_{\o
X}g$. The $\pi$-vector field $\o X$ is said to be  $\o d$-closed
{\em(resp. $\o d$-exact)} if its associated $\pi$-form $\omega$ is
$\o d$-closed {\em(resp. $\o d$-exact)}.
\end{defn}

The following result gives a simple characterization of $\o
d$-closed $\pi$-vector fields.\linebreak \vspace{-0.6cm}
\begin{thm}\label{2th.1}A $\pi$-vector field $\o X$ is $\o d$-closed if and only
if the operator $A_{\o X}$ is\linebreak self-adjoint.
\end{thm}

%\begin{proof}
\prof Let $\omega$ be the $\pi$-form associated to the $\pi$-vector
field $\o X$ under the duality defined by the Finsler metric $g$. By
(\ref{2eq.3}), we have\vspace{-0.2cm}
\begin{eqnarray*}
   (\o d \omega)(\o Y, \o Z)&=& \beta \o Y \cdot \omega(\o Z)-
    \beta \o Z \cdot \omega(\o Y)- \omega(\rho[\beta \o Y, \beta \o Z]) \\
   &=& \beta \o Y \cdot g(\o X, \o Z)-
    \beta \o Z \cdot g(\o X,\o Y)- g(\o X,\rho[\beta \o Y, \beta \o Z]) \\
   &=&  g(A_{\o X}\o Y, \o Z)+g(\o X,\nabla_{\beta \o Y}\o Z)-
   g(A_{\o X}\o Z,\o Y)-g(\o X,\nabla_{\beta \o Z}\o Y)\\
  & &- g(\o X,\rho[\beta \o Y, \beta \o Z])\\
   &=&g(A_{\o X}\o Y, \o Z)-g(A_{\o X}\o Z, \o Y)+g(\o X, Q(\o Y,\o
   Z)).
\end{eqnarray*}
As the horizontal torsion tensor $Q$ of the Cartan connection
$\nabla$ vanishes, the result follows. \ \ $\Box$
%\end{proof}

\begin{defn}The gradient of a function $f\in \mathfrak{F}(TM)$ is the $\pi$-vector
field $\o X$\linebreak defined by{\em\,:}\vspace{-0.2cm}
\begin{equation}\label{2eq.11}
    i_{\o X}\,g=\o d f=df o \beta . \vspace{-0.2cm}
\end{equation}
The gradient of the function $f$ is denoted by $grad\, f$.
\par A $\pi$-vector field $\o X$ is said to be  a gradient $\pi$-vector field
if it is the gradient of some function $f\in\mathfrak{F}(TM):\, \o X
= grad\, f$.
\end{defn}

In Riemannian geometry, it is well known that (\cite{r80},
\cite{r82}, \cite{r83}) the gradient of any function $f\in\f$ is a
closed vector field or, equivalently, $A_{ X}$ is self-adjoint. This
result is not in general true in Finsler geometry. This is again due
to the fact that $\o d^{2}\neq0$. Nevertheless, we
have\vspace{-0.2cm}

\begin{thm}\label{2th.2} The following assertions are equivalent{\,\em:}
\vspace{-0.2cm}
\begin{description}
    \item[(a)] The gradient of any  function $f\in\mathfrak{F}(TM)$ is  $\o
d$-closed, or, equivalently, gradient $\pi$-vector fields are $\o
d$-closed.
    \item[(b)]  The curvature tensor $\mathfrak{R}$ of the Barthel
    connection vanishes.
    \item[(c)] The horizontal distribution is completely integrable.
    \item[(d)] The $\pi$-exterior  derivative $\o d$ has the property that $\o
d^{2}=0$.
\end{description}
\end{thm}
%\begin{proof}
\prof\\ $\textbf{(a)}\Longleftrightarrow\textbf{(b)}:$ Let
$f\in\mathfrak{F}(TM)$ be any arbitrary function and let $\o
X:=grad\, f$. For all $\o Y, \o Z\in\cp$ and for all $\o X=grad\,
f$, we have\vspace{-0.2cm}
\begin{eqnarray*}
  g(A_{\o X}\o Y,\o Z) - g(A_{\o X}\o Z,\o Y)&=&
  g(\nabla_{\beta\o Y}\o X,\o Z)-g(\nabla_{\beta\o Z}\o X,\o Y)\\
   &=&\beta \o Y\cdot g(\o X,\o Z)- g(\o X,\nabla_{\beta\o Y}\o Z)\\
   & & -\beta \o Z\cdot g(\o X,\o Y)+g(\o X,\nabla_{\beta\o Z}\o Y) \\
   &=& \beta \o Y\cdot ((df o \beta) \o Z)-\beta \o Z\cdot ((df o \beta) \o Y)\\
   && -(df o \beta) (\nabla_{\beta\o Y}\o Z-\nabla_{\beta\o Z}\o Y) \\
   &=&\beta \o Y\cdot ( \beta \o Z\cdot f)-\beta \o Z\cdot ( \beta \o Y\cdot f)-
    \beta (\nabla_{\beta\o Y}\o Z-\nabla_{\beta\o Z}\o Y)\cdot f\\
    &=& ([\beta\o Y, \beta\o Z]-\beta (\nabla_{\beta\o Y}\o Z-
    \nabla_{\beta\o Z}\o Y))\cdot f\vspace{-0.2cm}
\end{eqnarray*}
In view of (\ref{2eq.5}), the last equation takes the
form\vspace{-0.2cm}
\begin{equation}\label{2eq.12}
g(A_{\o X}\o Y,\o Z) - g(A_{\o X}\o Z,\o Y)=\gamma (R(\o Y,\o Z)\o
\eta)\cdot f=\gamma(\widehat{R}(\o Y,\o Z))\cdot f \,\,\,\,
\,\forall f\in\mathfrak{F}(TM) \vspace{-0.2cm}
\end{equation}
\par
Now, the required equivalence  follows from (\ref{2eq.12}), taking
into account Theorem \ref{2th.1}, Lemma \ref{lem.bracket} and the
fact that $\gamma$ is a monomorphism.\\
$\textbf{(b)}\Longleftrightarrow\textbf{(c)}:$ This equivalence
follows immediately from the identity
\,\cite{r85}:\newline\vspace{-0.3cm}
$$\mathfrak{R}(X,Y)=-v[hX,hY].\vspace{-0.4cm}$$\\
$\textbf{(c)}\Longleftrightarrow\textbf{(d)}:$ For all $f\in
\mathfrak{F}(TM) $ and $\o X,\o Y \in\cp$, we have\vspace{-0.2cm}
$$({\o d}^{2}f)(\o X,\o Y)=\beta\o X\cdot(\beta \o Y \cdot f)-
\beta\o Y\cdot(\beta \o X \cdot f)-\beta(\rho[\beta \o X,\beta \o
Y])\cdot f=(I-\beta \, o\, \rho)\cdot [\beta \o X,\beta \o Y]\cdot
f.\vspace{-0.2cm}$$ Hence, we have $({\o d}^{2}f)(\o X,\o
Y)=v[\beta \o X,\beta \o Y]\cdot f \, \,\, \forall\,f\in\mathfrak{F}(TM)$. \\
The above equation, shows that the horizontal distribution is
completely integrable if and only if ${\o d}^{2}=0$ on
$\mathfrak{F}(TM)$. This proves the implication
$(d)\Longrightarrow(c)$. For the proof of the converse implication,
refer to \,\cite{r77}. \ \ $\Box$
%\end{proof}

\begin{defn}\label{2def.5} A Finsler manifold $(M,L)$ is said to be
of scaler curvature $\kappa$ if the $(v)h$-torsion tensor
$\widehat{R}$ is written in the form{\,\em:} \vspace{-0.2cm}
\begin{equation}\label{2eq.13}
   \widehat{R}:=R \otimes \o \eta =\omega \wedge \phi, \vspace{-0.2cm}
\end{equation}
where $\phi= I-L^{-1} \ell \otimes \o \eta$;\, $\ell:=L^{-1}i_{\o
\eta}\,g=dL\,o\,\gamma$ and $\omega:=\frac{1}{3}L(L\o
d_{J}\kappa+3\kappa\, \ell)$; \linebreak $\o d_{J}\kappa$ is the
$\pi$-lift of $d_{J}\kappa=d\kappa\, o\, J$. \em{(}$\kappa$ being
the horizontal scalar curvature).
\end{defn}

The formula (\ref{2eq.13}) is found in\ {\cite{r78}} in a local
coordinate from. It expresses the characteristic property of a
Finsler space of scalar curvature.

\begin{prop}Let $(M,L)$ be a Finsler manifold of scalar curvature $\kappa$.
The gradient of any  function $f\in\mathfrak{F}(TM) $ is $\o
d$-closed {\em({\it equivalently, all gradient $\pi$-vector fields
are $\o d$-closed}\,)} if  the horizontal scaler curvature $\kappa$
vanishes.
\end{prop}

%\begin{proof}
\prof The result follows from (\ref{2eq.12}) and (\ref{2eq.13}),
taking Theorem \ref{2th.1} into account. \ \ $\Box$
%\end{proof}

\begin{thm}\label{2th.3}Let $(M,L)$ be a Finsler manifold of nonzero scalar curvature
$\kappa$. Let $f\in\mathfrak{F}(TM) $ be an everywhere-nonzero
positively homogeneous function of degree $r$ in the directional
arguments. If $grad\,f$ is $\o d$-closed, then $f(x,y)=h(x)L^{r}$,
for some function $h\in\mathfrak{F}(M)$.
\end{thm}

%\begin{proof}
\prof Let $\o X:=grad\,f$. As $\o X$ is $\o d$-closed, then by
Theorem \ref{2th.1} and (\ref{2eq.12}), we have\vspace{-0.2cm}
$$\gamma \,(R(\o Y,\o Z)\o
\eta)\cdot f =0. \vspace{-0.2cm}$$ But since $(M,L)$ is of scaler
curvature, then \vspace{-0.2cm}
$$\gamma \,((\omega \wedge \phi )(\o Y,\o Z))\cdot f=0 . \vspace{-0.2cm}$$ Setting $\o Z=\o \eta$ in the
above equation, noting that $\phi(\o \eta)=0$, we get\vspace{-0.2cm}
$$ \omega(\o \eta) (\gamma \, \phi (\o Y))\cdot f=0 . \vspace{-0.2cm}$$
As $\omega(\o \eta)=\frac{1}{3}L^{2}(\mathcal{C}\cdot \kappa
+3\kappa)\neq0$, it follows that\vspace{-0.2cm}
$$ (\gamma \, \phi (\o Y))\cdot f=0 . \vspace{-0.2cm}$$
Consequently, by the definition of $\phi$, we get\vspace{-0.2cm}
\begin{equation}\label{2eq.14}
     \gamma \o Y\cdot f-L^{-1}\ell(\o Y)\,\mathcal{C}\cdot f=0 ,\vspace{-0.2cm}
\end{equation}
from which, since $f$ is homogeneous of degree $r$,\vspace{-0.3cm}
$$  df\, o \, \gamma -rf L^{-1} \ell= 0,\,\,
\text{or}\ \ \ \  \frac{df}{f}\, o \, \gamma -r\frac{dL}{L}\, o
\,\gamma = 0.\vspace{-0.4cm}$$ This equation is equivalent
to\vspace{-0.2cm}
$$ d(\log(fL^{-r}))\, o \,\gamma =0 .\vspace{-0.2cm}$$
Hence, $fL^{-r}$ is independent of the directional arguments and,
consequently, there exists a function $h\in\mathfrak{F}(M)$ so that
$f=h(x)L^{r}$. \ \ $\Box$
%\end{proof}

\begin{cor} Under the hypothesis of Theorem \ref{2th.3}, if the function $f\in\mathfrak{F}(TM)$
is homogeneous of degree zero in the directional arguments, then $f$
is a function of positional arguments only, that is $f\in\f.$
\end{cor}
\vspace{-5pt}
\par In fact, this follows from (\ref{2eq.14} ), since in this case
$\mathcal{C} \cdot f=0$.

\vspace{8pt}
\par
Let $U$ be an open subset of $\T M$. An assignment\vspace{-0.2cm}
$$\mathfrak{D}:u\in U\longrightarrow\mathfrak{D}_{u}\subset
P^{-1}(u)=\{u\}\times T_{\pi(u)}M,\vspace{-0.32cm}$$ such that every
$\mathfrak{D}_{u}$ is an $m$-dimensional vector subspace of
$P^{-1}(u)$, is called an\linebreak $m$-dimensional
$\pi$-distribution on $U$. If $\o X$ is a $\pi$-vector field on $U$
with $\o X(u)\in\mathfrak{D}_{u}$ for every $u\in U$, we say that
$\o X$ belongs to $\mathfrak{D}$ and we write $\o X\in
\mathfrak{D}$.
\par
For a given regular connection on $\p$, an $m$-dimensional
$\pi$-distribution $\mathfrak{D}$ on $U$ is said to be $h$-
involutive if for every $\pi$-vector fields $\o X$ and $\o Y$,
$\rho[\beta \o X, \beta \o Y]$ belongs to $\mathfrak{D}$ whenever
$\o X$ and $\o Y$ belong to $\mathfrak{D}$.

\begin{defn} Let $\o X$ be a given $\pi$-vector field on an
open subset $U$ of $\tm$. The $(n-1)$-dimensional $\pi$-distribution
\vspace{-0.2cm}
$$\mathfrak{D} : u\in U\longmapsto\mathfrak{D}_{u}:=\{\o Y\in P^{-1}(u): g(\o X,\o
Y)=0\vspace{-0.2cm}\}$$ is called the $\pi$-distribution generated
by (or associated with) $\o X$.
\end{defn}

\begin{thm}If $\o X\in \cp$ is $\o d$-closed, then the
$\pi$-distribution $\mathfrak{D}$ generated by $\o X$ is
$h$-involutive.
\end{thm}

%\begin{proof}
\prof Suppose that $\o Y$ and $\o Z$ belong to $\mathfrak{D}$. As
the $h$-torsion tensor of the Cartan connection $\nabla$ vanishes,
then\vspace{-0.2cm}
\begin{eqnarray*}
  g(\rho[\beta \o Y,\beta \o Z], \o X) &=& g(\nabla_{\beta \o Y}\o
  Z-\nabla_{\beta \o Z}\o Y,\o X) \\
  &=& \beta \o Y \cdot g(\o Z,\o X)-g(\o Z,\nabla_{\beta \o Y}\o X)-
  \beta \o Z \cdot g(\o Y,\o X)+ g(\o Y,\nabla_{\beta \o Z}\o X).\vspace{-0.2cm}
\end{eqnarray*}
Since $g(\o X,\o Y)=0=g(\o X,\o Z)$, it follows that\vspace{-0.2cm}
$$ g(\rho[\beta \o Y,\beta \o Z], \o X)=g(A_{\o X}\o Z,\o Y)-g(\o Z,A_{\o X}\o Y).
\vspace{-0.2cm}$$ Then, by Theorem \ref{2th.1}, $\mathfrak{D}$ is
$h$-involutive. \ \ $\Box$
%\end{proof}

\vspace{9pt}
\par
The following result characterizes certain $\o d$-closed
$\pi$-vector fields. \vspace{-0.1cm}
\begin{prop} Let $\o X$ be a $\pi$ vector field such that $i_{\o X}\,\o dg=0$.
Then, $\o X$ is $\o d$-closed if and only if $\beta\o X$ is an
isometry.
\end{prop}
\vspace{-8pt}
\par
The result follows directly from the identity $\mathfrak{L}_{\beta\o
X}\,g=i_{\o X}\,\o dg+\o di_{\o X}\,g$, where $\mathfrak{L}_{X}$ is
the Lie derivative with respect to $X\in \mathfrak{X}(TM)$
\,\cite{r77}. \vspace{9pt}
\par
A generalized Randers manifold \ \cite{r79} is a Finsler manifold
$(M,L^{*})$ whose Finsler structure $L^{*}$ is given by
$L^{*}=L+\alpha$, where $L$ is a Finsler structure on $M$ and
$\alpha=g(\o b,\o \eta)$; $\o b$ being the $\pi$-vector field
defined in terms of a given $1$-form $\delta$ on $M$
by\linebreak\vspace{-0.2cm}
$$ \delta(X)= g(\o b,\o X)\ \ \forall\,X\in \mathfrak{X}(M),\vspace{-0.2cm}$$
where $\o X$ is the $\pi$-life of $X$. \\The change
$L\longrightarrow L^{*}=L+\alpha$\, is called a generalized Randers
change.

\par
The next result gives a characterization of $\o d$-closeness of a
remarkable $\pi$-vector field $\o m$ associated with Randers
changes.
\begin{prop} Let $L^{*}=L+\alpha$ be a generalized Randers change
with closed $\alpha$. The $\pi$-vector field $\,\o m:=\o b
-L^{-2}\alpha\,\o \eta\,$ is $\o d$-closed in the Finsler manifold
$(M,L^{*})$ if and only if $\,\tau\o m\,$ is $\o d$-closed in the
Finsler manifold $(M,L)$, where $\tau=L^{*}L^{-1}$.
\end{prop}

%\begin{proof}
\prof  The relation between the Finsler metrics $g$ and $g^{*}$
associated with the Finsler structures $L$ and $L^{*}$ is given by \
\cite{r79}\,:\vspace{-0.2cm}
$$g^{*}=\tau(g-\ell \otimes \ell)+ \ell^{*} \otimes \ell^{*},\vspace{-0.2cm}$$
where $\tau=L^{*}L^{-1}, \ell= dL\, o\, \gamma$ and $\ell^{*}= \ell
+ d\alpha\, o\, \gamma$.\\
Now, if $d\alpha=0$, then $\ell^{*}=\ell$ and
consequently\vspace{-0.2cm}
$$g^{*}=\tau g+ (1-\tau)\ell \otimes \ell.\vspace{-0.28cm}$$
Then\vspace{-0.1cm}
$$i_{\o m}g^{*}=\tau i_{\o m}g+ (1-\tau)\ell(\o m) \ell.\vspace{-0.07cm}$$
As $\ell(\o m)=0$, as one can easily show, then the above relation
reduces to $\,i_{\o m}\,g^{*}= i_{\tau\o m}\,g$, from which the
result follows. \ \ $\Box$
%\end{proof}

\vspace{8pt}

\par
Two Finsler structure $L$ and $\tilde{L}$ on a manifold $M$ are said
to be conformal\ \cite{r62} if $\tilde{L}=e^{\sigma(x)}L$, for some
function $\sigma(x)$ on $M$. The transformation
$L\longrightarrow\tilde{L}=e^{\sigma(x)}L$ is called a conformal
transformation (or a conformal change). The relation between the
Finsler metrics $g$ and $\tilde{g}$ associated with $L$ and
$\tilde{L}$ respectively is given by\,: $\tilde{g}=e^{2\sigma(x)}g$.
\par
Let $\o X$ be a $\pi$-vector field, then $i_{\o X}\,\tilde{g}=
e^{2\sigma(x)}i_{\o X}\,g$, and so \vspace{-8pt}
$${\o d }i_{\o X}\,\tilde{g}=
e^{2\sigma(x)}{\o d }i_{\o X}\,g+2e^{2\sigma(x)}(d \,o \,
\beta)\sigma \otimes i_{\o X}\,g=e^{2\sigma(x)}{\o d }i_{\o
X}\,g+2e^{2\sigma(x)}\frac{\partial\sigma}{\partial x^{k}} \,dx^{k}
\otimes i_{\o X}\,g.\vspace{-0.2cm}$$ Consequently, $\,{\o d }i_{\o
X}\,\tilde{g}= e^{2\sigma(x)}{\o d }i_{\o X}\,g\,$ if and only if
$\,\frac{\partial\sigma}{\partial x^{k}}=0$, or equivalently, if and
only if $\sigma(x)$ is a constant function (provided that $M$ is
connected). This proves the following result, which gives a
characterization of homotheties in terms of $\o d$-closed
$\pi$-vector fields.

\begin{thm} A $\o d$-closed $\pi$-vector field remain $\o d$-closed
under a conformal transformation if and only if this transformation
is a homothety.
\end{thm}

\providecommand{\bysame}{\leavevmode\hbox
to3em{\hrulefill}\thinspace}
\providecommand{\MR}{\relax\ifhmode\unskip\space\fi MR }
\providecommand{\MRhref}[2]{%
  \href{http://www.ams.org/mathscinet-getitem?mr=#1}{#2}
} \providecommand{\href}[2]{#2}

\end{document}